\documentclass[preprint,12pt]{elsarticle}

\newif\ifarxiv
\arxivtrue

\usepackage{amsmath,amssymb,amsfonts}
\usepackage{booktabs}
\usepackage{graphicx}
\usepackage{longtable}
\usepackage{xurl}
\ifarxiv
  \usepackage[colorlinks,linkcolor=blue,citecolor=blue,urlcolor=blue]{hyperref}
\else
  \usepackage[hidelinks]{hyperref}
\fi

\biboptions{sort&compress}

\newcommand{\figw}{0.82\linewidth}

\journal{Journal of Energy Storage}

\begin{document}

\begin{frontmatter}

\title{Co-Optimized Generation, Transmission, and Storage Expansion: System Value and Optimal Duration of Pumped-Storage Hydropower}

\author[psr]{Rafael Benchimol Klausner}
\author[psr]{Rafael Kelman}

\affiliation[psr]{organization={PSR},
  addressline={Praia de Botafogo 370, Botafogo},
  city={Rio de Janeiro},
  state={RJ},
  postcode={22250-040},
  country={Brazil. Corresponding author: Rafael Benchimol Klausner,
           rafabench@psr-inc.com}}

\begin{abstract}
Expansion planning models usually fix storage duration before optimization,
setting how much storage power to build but not for how long it can
discharge. We present a generation--transmission--storage expansion framework
in which candidates of many durations compete on annualized cost, making the
duration mix an optimization outcome. It is a rolling-horizon, two-stage
stochastic linear program, each five-year stage is a full year at hourly
resolution under ten coherent inflow, wind, and solar scenarios. We apply it
to the Brazilian Interconnected System over 2030--2050, where demand roughly
doubles to 1{,}716~TWh/year and variable renewable energy (VRE) supplies most
new capacity. Two cases are compared: 36 pumped-storage hydropower (PSH)
candidates over four subsystems and nine durations (4--144~h) alongside 4-h
battery energy storage systems (BESS), and BESS alone.

With PSH available, the model builds 31.2~GW / 755~GWh of PSH and no BESS,
dominated by 12-h capacity with 4.9~GW of 72-h units. Without PSH it builds
38~GW / 152~GWh of BESS and 7.8~GW more gas-fired capacity, mostly
open-cycle peakers. PSH lowers 2050 annualized system cost by
US\$~5.0~billion/year (6.8\%), operating cost by 23.5\%, thermal generation
by 34~TWh/year, and long-run marginal cost by 20\%; 2035 VRE curtailment
falls from 8.4\% to 3.2\%. The magnitudes are specific to Brazil, but the underlying mechanism is general: 
limiting storage candidates to a single duration understates the system’s optimal energy-storage requirement
 and overstates its residual need for thermal capacity.
\end{abstract}



\end{frontmatter}

\section*{Nomenclature}
\addcontentsline{toc}{section}{Nomenclature}

\begin{longtable}{@{}p{0.20\linewidth}p{0.76\linewidth}@{}}
\multicolumn{2}{@{}l@{}}{\textit{Sets and indices}}\\
\\
$t \in \mathcal{T}$ & Hours, $\mathcal{T}=\{1,\dots,8760\}$\\
$\omega \in \Omega$ & Operating scenarios, with probabilities $p_\omega$\\
$s \in \mathcal{S}$ & Subsystems (electrical zones)\\
$f \in \mathcal{F}$ & Renewable plants; $\mathcal{F}_s$ those in subsystem $s$\\
$i \in \mathcal{I}$ & Thermal plants; $\mathcal{I}_s$ those in subsystem $s$\\
$h \in \mathcal{H}$ & Aggregated (equivalent) hydro reservoirs\\
$b \in \mathcal{B}$ & Storage units (BESS and PSH)\\
$l \in \mathcal{L}$ & Interconnections; $\delta^{\mathrm{in}}(s)$,
$\delta^{\mathrm{out}}(s)$ those entering and leaving $s$\\
$p \in \mathcal{P}$ & Candidate expansion projects; $\mathcal{P}(\cdot)$
those attached to a given plant, storage unit or corridor\\
$y$ & Planning year of the rolling horizon\\
$m(t)$ & Month containing hour $t$\\[4pt]
\\
\multicolumn{2}{@{}l@{}}
{\textit{Parameters}}\\
\\
$A_p$ & Annualized fixed cost of project $p$ (US\$/MW/yr)\\
$\mathrm{CRF}_p$ & Capital recovery factor, from discount rate $r_p$ and
economic lifetime $n_p$\\
$c^{\mathrm{T}}_i$ & Thermal variable operating cost (US\$/MWh)\\
$c^{\mathrm{B}}_b$ & Storage O\&M throughput cost (US\$/MWh)\\
$c^{\mathrm{C}}_f$ & Renewable curtailment penalty (US\$/MWh)\\
$c^{\mathrm{D}}$ & Cost of energy not served (US\$/MWh)\\
$D_{s,t}$ & Demand of subsystem $s$ in hour $t$ (MW)\\
$\rho_{f,t,\omega}$ & Hourly capacity factor of renewable plant $f$\\
$\bar{G}_f$ & Existing capacity of renewable plant $f$ (MW)\\
$\bar{G}^{\mathrm{av}}_{i,t}$ & Available (derated) thermal capacity (MW)\\
$\underline{\gamma}_{i,m(t)}$ & Monthly must-run share of thermal plant $i$\\
$a_{h,t,\omega}$ & Scenario inflow to reservoir $h$, in energy (MWh)\\
$\underline{E}_h$, $\overline{E}_h$ & Reservoir storage bounds (MWh)\\
$\underline{o}_h$ & Minimum environmental outflow (MW)\\
$\eta^{+}_b$, $\eta^{-}_b$ & Charge and discharge efficiencies\\
$\tau_b$ & Storage duration, i.e.\ energy-to-power ratio (h)\\
$\delta_b$ & Depth-of-discharge reserve\\
$\eta^{\mathrm{L}}$ & Interconnection loss factor\\
$\overline{F}^{\mathrm{from}}_l$, $\overline{F}^{\mathrm{to}}_l$ & Existing
transfer limits per direction (MW)\\
$\overline{X}_p$ & Resource potential of project $p$ (MW)\\[4pt]
\\
\multicolumn{2}{@{}l@{}}{\textit{Decision variables}}\\
\\
$x_p$ & First stage: capacity built for project $p$ (MW)\\
$X_b$, $X_l$ & Invested capacity aggregated at storage unit $b$ /
corridor $l$ (MW)\\
$g^{\mathrm{R}}_{f,t,\omega}$, $c_{f,t,\omega}$ & Renewable dispatch and
curtailment (MW)\\
$g^{\mathrm{T}}_{i,t,\omega}$ & Thermal dispatch (MW)\\
$g^{\mathrm{H}}_{h,t,\omega}$, $v_{h,t,\omega}$ & Hydro generation and
spillage (MW)\\
$E_{h,t,\omega}$ & Hydro reservoir storage (MWh)\\
$q^{+}_{b,t,\omega}$, $q^{-}_{b,t,\omega}$ & Storage charge and discharge
(MW)\\
$e_{b,t,\omega}$ & Storage state of charge (MWh)\\
$f_{l,t,\omega}$ & Interconnection flow (MW)\\
$r_{s,t,\omega}$ & Load shedding (MW)\\
\end{longtable}

\section{Introduction}
\label{sec:intro}

High shares of variable renewable energy (VRE) split the balancing problem
across different time scales. Solar ramps must be covered within the day, wind droughts
persist for several days, and a seasonal imbalance appears wherever hydro
inflows, wind, or demand follow an annual cycle. Storage can serve all three,
but a
fleet sized for one is not sized for the others, because each service calls
for a different energy-to-power ratio. The planning question is therefore not
only how much storage power to build, but which durations should compete and
how the selected mix changes generation, transmission, and system operation.
Expansion models commonly settle that question by assumption, admitting
storage as a single class whose duration is fixed before
optimization~\cite{bistline2020}.

Pumped-storage hydropower (PSH) is the most mature long-duration storage
technology and accounts for most installed storage capacity
worldwide~\cite{iha2024,blakers2021}. Variable-speed and ternary
configurations have expanded its operating
range~\cite{koritarov2022,zhao2024,papadakis2023}, and global siting surveys
report closed-loop resource potential far in excess of any plausible
system requirement~\cite{stocks2021}. Planning generation, transmission, and
storage in separate steps can undervalue PSH because it provides energy
arbitrage, capacity, curtailment reduction, and transmission relief
jointly~\cite{krishnan2016,jacobson2023}.

This paper presents a co-optimized
generation--\allowbreak transmission--\allowbreak storage (G,T\&S) expansion
framework in which storage candidates spanning a wide range of durations
compete against one another, so the duration mix of the selected portfolio is
left to the optimization rather than assumed in advance. Every candidate
class is governed by the same hourly state-of-charge equations and differs
only in its energy-to-power ratio $\tau_b$, efficiency, usable depth of
discharge, annualized cost, and deployment limit. Because capacity is chosen
per class, the optimizer sets storage power, energy ($\tau_b X_b$), and
location jointly with generation and transmission investment and with hourly
dispatch. The model is a
rolling-horizon, two-stage stochastic linear program, implemented in
Julia/JuMP~\cite{jump} and solved with FICO Xpress~\cite{xpress}. Nothing in
the formulation is specific to a particular system.

We exercise the framework on the Brazilian Interconnected System (SIN) over
2030--2050, which is a demanding case for duration-differentiated storage on
four counts. Its generation mix is already largely renewable, and wind and
solar growth has outpaced the expansion of the grid and of the flexibility
resources needed to absorb it~\cite{iea_brazil2025}. Firm hydropower
additions are constrained by socio-environmental limits, so VRE supplies most
of the growth to 1{,}716~TWh/year by 2050 (Section~\ref{sec:case}). The
storage capability of its reservoirs has declined relative to load. And two
combined-cycle contract families carry must-run obligations in opposite
halves of the year, so inflexible thermal generation is present year-round.
PSH has entered Brazilian planning studies as a candidate technology, but the
official ten-year expansion plan has not selected any PSH
capacity~\cite{epe_pde2034}, and the first capacity-reserve auction for
storage was opened to batteries alone~\cite{mme136}.

We compute two expansion plans with identical demand, scenarios, renewable
and thermal candidates, and transmission reinforcements. Only the storage
candidate set differs. In the \emph{with-PSH} case, 36 PSH candidates
comprising nine storage durations (4--144~h) in each of the four SIN
subsystems compete with 4-h lithium-ion battery energy storage systems
(BESS). In the \emph{no-PSH} case, only the BESS candidates are available.

The case study addresses three storage-planning questions:
\begin{enumerate}
  \item Which PSH power, energy, and duration mix is selected when 4--144-h
  PSH candidates compete with 4-h BESS?
  \item How does admitting PSH change the optimal storage, generation, and
  transmission portfolios and their annualized system cost?
  \item Which operating patterns, from the daily cycle to multi-day and
  seasonal events, explain the resulting changes in curtailment, thermal
  generation, adequacy, and marginal cost?
\end{enumerate}

Because the two cases use identical non-storage assumptions, their
differences isolate the effect of admitting duration-differentiated PSH
candidates, subject to the stated technology costs and deployment limits.

\section{Related work}
\label{sec:related}

\subsection{Storage in capacity expansion models}

Capacity expansion models select generation, transmission, and storage
investments by minimizing investment plus expected operating cost. Reviews of
co-optimization practice describe how resolving those decisions sequentially
rather than jointly changes the resulting portfolio, and catalogue the
modeling approaches available for treating them together~\cite{krishnan2016}.
Placing storage in such a model requires chronological
operation and explicit state-of-charge dynamics: an exogenous capacity credit
applied to a load-duration curve cannot express how a fleet's contribution
changes with its energy-to-power ratio~\cite{bistline2020}. The level of
temporal and operational detail is itself a modeling choice with quantitative
consequences, since coarse time slices and relaxed operating constraints
misrepresent the flexibility a system can actually
deliver~\cite{poncelet2016,helisto2019}. Open-source platforms such as PyPSA
and Switch have made hourly, chronological expansion planning
practicable at national scale~\cite{pypsa2018,switch2019}.

Studies of long-duration energy storage (LDES) in deep-decarbonization
settings show that storage value is highly sensitive to the energy-to-power
ratio and that duration classes serve qualitatively different system
functions, from intra-day shifting to seasonal smoothing. For durations
beyond a day, the cost per kWh of the storage medium becomes a decisive
parameter~\cite{sepulveda2021,denholm2021}. Dowling
et al.~\cite{dowling2020} report that system cost responds roughly twice as
strongly to the cost of long-duration storage as to the cost of batteries in
wind--solar systems, which makes the duration composition of the candidate
set a first-order modeling assumption.
Reporting storage only by power capacity can obscure the energy
volume and cycling timescale that produce its system value. Letting several
duration classes compete allows those attributes to be evaluated within the
capacity mix rather than imposed before optimization.
Two-stage and multi-stage stochastic programming is the standard framework
for capturing operational uncertainty in these
decisions~\cite{birge2011,pereira1991}.

\subsection{PSH technology, resource, and valuation}

PSH accounts for most installed storage power and energy capacity
worldwide~\cite{iha2024,blakers2021}, and reviews of its development across
the major electricity markets describe how market structure and remuneration
arrangements have shaped where it was built~\cite{barbour2016}. Recent technology reviews examine variable-speed,
ternary, and quaternary configurations, which extend PSH services to fast
frequency response and pumping-power
modulation~\cite{koritarov2022,zhao2024,papadakis2023,dong2020,ieee_tr134}.

Resource assessments support that reading. Global GIS screening identifies
about 616{,}000 closed-loop, off-river sites with roughly 23{,}000~TWh of
combined storage potential, orders of magnitude beyond any plausible
system requirement~\cite{stocks2021}, and seasonal configurations that
exploit existing river topography extend the accessible timescale to
months~\cite{hunt2020}. Bottom-up cost models show an economy of duration:
reservoir costs scale with storage volume while powerhouse costs are largely
fixed, so the unit cost of stored energy falls with
duration~\cite{cohen2023,jacobson2023}.

Where PSH has been placed inside an expansion or dispatch model, the value it
returns depends on how its candidates are represented. Li
et al.~\cite{li2023} quantify the curtailment and emissions reductions
attainable with PSH in northwest China, and Zhu et al.~\cite{zhu2026} find
for the Chinese system that aggregated representations of closed-loop PSH
overestimate the capacity required by about 30\% relative to a
reservoir-level formulation. Valuation studies also show that
energy-arbitrage revenues alone can understate PSH system value, and that
integrated G,T\&S planning can quantify capacity, ancillary-service,
curtailment-avoidance, and transmission-relief benefits
jointly~\cite{jacobson2023}. The relevant comparison is therefore not
limited to storage cost or arbitrage revenue; it also includes changes in
generation, transmission, and system operation when storage technologies and
durations compete.

\subsection{Case-study context: Brazil}

Brazilian ten-year expansion planning has begun to incorporate storage, and
PDE~2034 includes storage technologies in its reference scenario for the
first time~\cite{epe_pde2034}. No PSH capacity has been selected in the
plan's optimized results, however, and the first capacity-reserve auction for
storage was designed for batteries only~\cite{mme136}. The duration question
is therefore live in Brazilian procurement before it has been settled in
Brazilian planning. Silva et al.~\cite{silva2025} compare PSH with
gas-fired peaking capacity in the Brazilian market and find PSH the cheaper
option for daily operating windows longer than about seven hours, a threshold
that bears directly on which duration classes an expansion model will
select.

Internationally, Australia's Integrated System Plan treats storage and transmission as substitutable options in a single
co-optimization, with long-duration storage backed by state-level
contract-for-difference schemes~\cite{aemo_isp}. China, the world's largest
PSH market, remunerates PSH through a two-part capacity and energy
tariff~\cite{iha2024}. For Brazil, Weber
et al.~\cite{weber2024} compared national expansion pathways with
and without seasonal PSH using a monthly-resolution MESSAGEix model under
climate-change scenarios. Compared with that work, this study co-optimizes
generation, transmission, and storage at hourly resolution with coherent
hydro--wind--solar scenarios and a duration-differentiated PSH candidate
set. The model determines the storage \emph{duration} mix and its interaction
with daily solar cycles and thermal must-run windows.

\section{Methodology}
\label{sec:method}

The formulation below is stated for a generic multi-zone hydro-thermal
system; Section~\ref{sec:case} instantiates it for the SIN. It compares
storage alternatives at the system level: PSH and
BESS use the same hourly storage equations and compete with generation and
transmission investments, while their candidate classes differ in duration,
efficiency, depth-of-discharge reserve, annualized cost, and deployment
limits. The optimizer selects storage power, energy
($\tau_b X_b$), location, and dispatch jointly.

The expansion planning model is a two-stage stochastic linear program.
First-stage variables are capacity investments; second-stage variables
describe hourly system operation in each scenario. The model represents
subsystems (electrical zones), aggregated hydro plants, thermal plants,
renewable plants, storage units (BESS and PSH), and interconnections. The Nomenclature lists the sets,
parameters, and variables used below.

\subsection{Objective}

The objective minimizes annualized investment cost plus expected operating
cost:
\begin{align}
\min\;
& \sum_{p\in\mathcal{P}} A_p\, x_p
  + \sum_{\omega\in\Omega} p_\omega
  \Bigg[
    \sum_{t,i} c^{\mathrm{T}}_i\, g^{\mathrm{T}}_{i,t,\omega}
    \nonumber\\
& \quad
    + \sum_{t,b} c^{\mathrm{B}}_b \big(q^{+}_{b,t,\omega}+q^{-}_{b,t,\omega}\big)
    + \sum_{t,f} c^{\mathrm{C}}_{f}\, c_{f,t,\omega}
    \nonumber\\
& \quad
    + \sum_{t,s} c^{\mathrm{D}}\, r_{s,t,\omega}
  \Bigg],
\label{eq:objective}
\end{align}
where $A_p = (\mathrm{CAPEX}_p \cdot \mathrm{CRF}_p +
\mathrm{OM}^{\mathrm{fix}}_p)$ is the annualized fixed cost of project $p$,
with capital recovery factor
$\mathrm{CRF}_p = r_p (1+r_p)^{n_p} / \big((1+r_p)^{n_p}-1\big)$ for
discount rate $r_p$ and economic lifetime $n_p$; $c^{\mathrm{T}}_i$ is the
thermal variable cost, $c^{\mathrm{B}}_b$ a storage O\&M throughput cost,
$c^{\mathrm{C}}_f$ an (optional) curtailment penalty, and $c^{\mathrm{D}}$
the cost of energy not served (1{,}666~US\$/MWh in this study).

\subsection{Constraints}

\subsubsection{Energy balance}
For every subsystem $s$, hour $t$, and scenario $\omega$:
\begin{align}
&\sum_{f\in\mathcal{F}_s} g^{\mathrm{R}}_{f,t,\omega}
+ \sum_{h\in\mathcal{H}_s} g^{\mathrm{H}}_{h,t,\omega}
+ \sum_{i\in\mathcal{I}_s} g^{\mathrm{T}}_{i,t,\omega}
+ \sum_{b\in\mathcal{B}_s}\big(q^{-}_{b,t,\omega}-q^{+}_{b,t,\omega}\big)
\nonumber\\
&\quad
+ \eta^{\mathrm{L}}\!\!\sum_{l\in\delta^{\mathrm{in}}(s)}\!\! f_{l,t,\omega}
- \!\!\sum_{l\in\delta^{\mathrm{out}}(s)}\!\! f_{l,t,\omega}
+ r_{s,t,\omega}
= D_{s,t},
\label{eq:balance}
\end{align}
with the interconnection loss factor $\eta^{\mathrm{L}}$ applied at the
receiving end.

\subsubsection{Renewables}
Available generation is split between dispatch and curtailment,
\begin{equation}
g^{\mathrm{R}}_{f,t,\omega} + c_{f,t,\omega}
 = \rho_{f,t,\omega}\,\big(\bar{G}_f + \textstyle\sum_{p\in\mathcal{P}(f)} x_p\big),
\label{eq:renewable}
\end{equation}
where $\rho_{f,t,\omega}$ is the hourly capacity factor profile and
$\bar{G}_f$ the existing capacity.

\subsubsection{Thermal plants}
Thermal dispatch is bounded by available capacity and by a monthly must-run
(inflexibility) share $\underline{\gamma}_{i,m(t)}$ representing
take-or-pay fuel contracts,
\begin{equation}
\underline{\gamma}_{i,m(t)}\,\bar{G}^{\mathrm{av}}_{i,t}
\;\le\; g^{\mathrm{T}}_{i,t,\omega} \;\le\; \bar{G}^{\mathrm{av}}_{i,t},
\label{eq:thermal}
\end{equation}
with $\bar{G}^{\mathrm{av}}_{i,t}$ the derated installed-plus-candidate
capacity.

\subsubsection{Hydro}
Aggregated reservoirs follow an energy balance with cyclic closure over the
year,
\begin{gather}
E_{h,t,\omega} = E_{h,t-1,\omega} + a_{h,t,\omega}
 - g^{\mathrm{H}}_{h,t,\omega} - v_{h,t,\omega},
\label{eq:hydro}\\
\underline{E}_h \le E_{h,t,\omega} \le \overline{E}_h, \qquad
E_{h,|\mathcal{T}|,\omega} = E_{h,0},\\
g^{\mathrm{H}}_{h,t,\omega} + v_{h,t,\omega} \ge \underline{o}_h ,
\end{gather}
where $a_{h,t,\omega}$ is the scenario inflow (in energy) and
$\underline{o}_h$ a minimum environmental outflow.

\subsubsection{Storage (BESS and PSH)}
Both technologies share one set of linear storage dynamics,
\begin{gather}
e_{b,t,\omega} = e_{b,t-1,\omega}
 + \eta^{+}_b\, q^{+}_{b,t,\omega}
 - \frac{1}{\eta^{-}_b}\, q^{-}_{b,t,\omega},
\label{eq:soc}\\
\delta_b\, \tau_b X_b \;\le\; e_{b,t,\omega} \;\le\; (1-\delta_b)\,\tau_b X_b,
\qquad
q^{\pm}_{b,t,\omega} \le X_b,
\label{eq:socbounds}\\
e_{b,|\mathcal{T}|,\omega} = e_{b,1,\omega},
\label{eq:soccyclic}
\end{gather}
where $X_b = \sum_{p\in\mathcal{P}(b)} x_p$ is the invested power capacity,
$\tau_b$ the duration (energy-to-power ratio, h), $\eta^{+}_b,\eta^{-}_b$
the charge/discharge efficiencies, and $\delta_b$ a depth-of-discharge
reserve. The cyclic condition~\eqref{eq:soccyclic} prevents fictitious net
energy drawdown over the year. The nine PSH candidate classes differ in the
duration parameter $\tau_b$. No operating service is assigned to a duration
class in advance; daily and multi-day roles emerge from the optimized hourly
dispatch.

\subsubsection{Transmission}
Interconnection flows are bounded by existing plus candidate capacity per
direction,
\begin{equation}
-\big(\overline{F}^{\mathrm{to}}_l + X_l\big) \;\le\; f_{l,t,\omega}
 \;\le\; \overline{F}^{\mathrm{from}}_l + X_l .
\label{eq:flow}
\end{equation}

\subsubsection{Investment bounds and rolling horizon}
Each project is bounded by its resource potential, $x_p \le \overline{X}_p$.
The 2030--2050 study is solved as a rolling horizon over the five
planning years $y \in \{2030,\,2035,\,\ldots,\,2050\}$. Year $y$ is solved as a complete
two-stage problem with year-specific demand, costs, and profiles, and its
optimal investments become lower bounds in the next planning year,
\begin{equation}
x_p^{(y+1)} \;\ge\; \sum_{y' \le y} \hat{x}_p^{(y')},
\label{eq:rh}
\end{equation}
so installed capacity is monotonically non-decreasing (no endogenous
retirements) and each five-year step re-optimizes both the incremental
build-out and the full hourly operation.

\subsection{Scenario generation}

Operational uncertainty enters through joint hydro--wind--solar scenarios.
Approximately 90 years of naturalized inflow records are processed by an
SDDP-based hydrological model~\cite{pereira1991} into 200 synthetic monthly inflow scenarios, which are reduced to ten representative scenarios by variance-preserving clustering. Hourly wind and solar capacity-factor profiles are synthesized from MERRA-2/ERA5 reanalysis data at 127 candidate measurement points, keeping hydro, wind, and solar realizations
chronologically coherent within each scenario so that inter-regional and
inter-variable correlations (e.g., dry years with strong Northeast winds)
are preserved. Each scenario is a full 8{,}760-hour year; the
scenario set preserves the joint distribution of solar ramps, multi-day
wind lulls, and wet/dry hydrology.

\subsection{Solution method}

The two-stage problem is solved as its \emph{deterministic equivalent}
(DE): all ten scenarios are assembled into a single LP with shared
first-stage variables and probability-weighted operating costs. The DE is
solved without a decomposition approximation and is tractable here because
ten hourly scenarios fit a single 128-GB compute node
(Section~\ref{sec:results:comp}). Each planning year's DE is solved with
the FICO Xpress barrier method with crossover disabled~\cite{xpress}.

\section{Case study: the Brazilian Interconnected System, 2030--2050}
\label{sec:case}

\begin{table}[!t]
\caption{Case study dimensions.}
\label{tab:dims}
\centering
\small
\begin{tabular}{@{}ll@{}}
\toprule
Subsystems & SE/CW, S, NE, N + Imperatriz hub \\
Hydro & 4 equivalent reservoirs, 112.2\,GW, 563\,TWh \\
Thermal plants & 296 existing/contracted \\
Renewable agents & 391 (84 expandable, $\sim$145\,GW potential) \\
Storage candidates & 4 BESS (4\,h) + 36 PSH (4--144\,h) \\
Transmission & 6 existing corridors + 4 candidates \\
Candidate projects & 221 (incl.\ thermal re-contracting) \\
Demand 2050 & 195.9\,GW avg, 243.0\,GW peak, 1{,}716\,TWh \\
Resolution & 8{,}760\,h/year, 10 scenarios \\
Planning years & 2030, 2035, 2040, 2045, 2050 \\
\bottomrule
\end{tabular}
\end{table}

The representation contains the four subsystems: Southeast/Center-West
(SE), South (S), Northeast (NE), and North (N). It also includes the
Imperatriz interconnection hub (NI) (Table~\ref{tab:dims}). Hydropower is
aggregated into
four equivalent reservoirs totaling 112.2~GW of capacity and 563~TWh
of useful storage, of which the Southeast holds roughly two thirds. The
thermal fleet comprises 296 plants, including the Angra~1--2 nuclear units.
Existing and contracted renewables amount to about 151~GW across 391 agents
(onshore/offshore wind, solar PV, small hydro, biomass, and distributed
generation). Demand follows the official ten-year plan (PDE~2034) through
2034~\cite{epe_pde2034} and is extended to 2050, reaching
1{,}716~TWh/year (196~GW average), 55\% of it in the Southeast.

Candidate thermal units include open-cycle gas peakers and two combined-cycle families with opposite seasonal must-run
regimes (pre-salt gas contracts, must-run December--May; LNG contracts,
must-run June--November), so inflexible thermal generation is present
year-round. This structure favors resources that can shift energy over
longer windows. Renewable curtailment carries a 30~US\$/MWh penalty in the
objective for all renewable plants, and energy not served is priced at the
value of lost load (VOLL) of 1{,}666~US\$/MWh; both appear inside reported
operating costs. All monetary values in this paper are stated in US dollars,
converted from Brazilian reais at 5~R\$/US\$.

\begin{table}[!t]
\caption{Demand projection by subsystem (2050).}
\label{tab:demand}
\centering
\small
\begin{tabular}{@{}lrrr@{}}
\toprule
Subsystem & Average (GW) & Peak (GW) & Energy (TWh/yr) \\
\midrule
Southeast/CW & 107.7 & 136.1 & 943.8 \\
South        &  36.0 &  53.9 & 315.2 \\
Northeast    &  34.7 &  41.6 & 303.8 \\
North        &  17.5 &  21.4 & 153.5 \\
\midrule
System total & 195.9 & 243.0 & 1{,}716.3 \\
\bottomrule
\end{tabular}
\end{table}

Table~\ref{tab:demand} details the 2050 demand projection by subsystem.
The transmission topology comprises six existing corridors (SE--S, SE--NI,
NE--NI, N--NI, SE--NE, and SE--N, with direction-dependent capacities of
3.0--9.7~GW, except N--NI, which is left unconstrained) and four candidate
reinforcements, each expandable continuously at its own annualized cost:
SE--S at 75.4, SE--NE at 79.0, SE--Imperatriz at 79.6, and NE--S at
85.8~US\$/kW/yr, over a 25-year economic lifetime. Inter-subsystem exchanges
incur loss factors. Because the Northeast concentrates the best wind resources while the
Southeast concentrates 55\% of demand, the SE--NE corridor and the storage
siting decision interact strongly, so the model co-optimizes storage siting
and SE--NE corridor expansion.

Storage candidates comprise one 4-h BESS per subsystem (83\% round-trip
efficiency, $\delta_b=0.20$) and, in the with-PSH case, 36 closed-loop
PSH candidates: nine duration classes (4, 6, 8, 12, 24, 48, 72, 100, and
144~h) in each of the four subsystems, with 80\% round-trip efficiency and
full depth of discharge, corresponding to minimum operative water levels. All storage costs enter the model as annualized fixed values (Table~\ref{tab:pshcost}). The PSH figures are site-level estimates produced with HERA-S, a prospecting framework that screens candidate sites and costs them from optimized dam-fill and excavation mass balances~\cite{albuquerque2026}. Annualized PSH cost grows sub-linearly with duration because upper-reservoir civil-works costs scale with storage volume while powerhouse costs are fixed~\cite{cohen2023}. A 144-h plant costs only 76\% more per kW than a 4-h plant but provides $36\times$ the
energy. PSH costs are held constant over the horizon, while
the BESS cost is assumed to decline linearly from 193.0~US\$/kW/yr in 2030 to
131.2~US\$/kW/yr in 2050, reflecting projected lithium-ion learning curves. By
2050, the 4-h BESS is 22\% cheaper per kW-year than the 4-h PSH.

\begin{table}[!t]
\caption{Annualized cost of storage candidates (2050 reference).}
\label{tab:pshcost}
\centering
\small
\begin{tabular}{@{}lcc@{}}
\toprule
Candidate & Duration (h) & Cost (US\$/kW/yr) \\
\midrule
BESS (Li-ion) & 4 & 131.2 \\
PSH & 4 & 168.2 \\
PSH & 6 & 176.4 \\
PSH & 8 & 180.4 \\
PSH & 12 & 184.6 \\
PSH & 24 & 208.6 \\
PSH & 48 & 228.4 \\
PSH & 72 & 234.8 \\
PSH & 100 & 260.2 \\
PSH & 144 & 295.6 \\
\bottomrule
\end{tabular}
\end{table}

Both cases include the same onshore and offshore wind, solar PV, small
hydro, biomass, gas thermal, thermal re-contracting, demand response, and
interconnection candidates. They differ \emph{only} in the storage
candidate set: with-PSH (36 PSH and 4 BESS candidates) versus no-PSH
(4 BESS candidates only).

\section{Results}
\label{sec:results}

\subsection{Expansion plans}
\label{sec:results:plans}

\begin{table}[!t]
\caption{Cumulative installed new capacity by technology (GW). Gas
candidates are shown separately for the open-cycle peakers and the two
combined-cycle contract families, whose must-run windows are opposite (LNG,
June--November; pre-salt, December--May).}
\label{tab:capacity}
\centering
\small
\begin{tabular}{@{}lrrrr@{\hskip 10pt}rrrr@{}}
\toprule
& \multicolumn{4}{c}{With PSH} & \multicolumn{4}{c}{No PSH} \\
\cmidrule(r{10pt}){2-5}\cmidrule{6-9}
Technology & '35 & '40 & '45 & '50 & '35 & '40 & '45 & '50 \\
\midrule
Onshore wind   & 20.7 & 58.3 & 91.3 & 122.4 & 23.1 & 58.7 & 91.7 & 122.7 \\
Solar PV       &      &  0.3 & 23.8 &  61.5 &      &      &  7.9 &  45.9 \\
Offshore wind  &      &      &      &   5.9 &      &      &  1.2 &   5.9 \\
PSH            &  8.4 & 10.7 & 18.2 &  31.2 &  --- &  --- &  --- &   --- \\
BESS (4\,h)    &      &      &      &       &  2.7 &  4.3 &  8.2 &  38.0 \\
Gas, open cycle &      &  2.4 &  2.4 &   2.4 &  3.1 &  3.6 &  5.5 &   7.8 \\
Gas CC, LNG      &      &  1.0 & 11.6 &  16.6 &      &  2.6 & 14.0 &  19.0 \\
Gas CC, pre-salt &      &      & 10.0 &  10.0 &      &      & 10.0 &  10.0 \\
Thermal re-contracting & 0.4 & 3.3 & 6.5 & 6.5 & 0.4 & 4.3 & 9.2 & 9.2 \\
Demand response &  2.2 &  2.8 &  3.2 &   3.6 &  2.2 &  2.8 &  3.2 &   3.6 \\
Small hydro    &      &  0.3 &  1.8 &   3.3 &      &  0.3 &  1.8 &   3.3 \\
Transmission   &  2.0 & 11.3 & 19.3 &  19.3 &  1.4 & 10.7 & 19.1 &  19.1 \\
\bottomrule
\end{tabular}
\end{table}

Admitting duration-differentiated PSH changes both the selected storage fleet
and the wider expansion plan. Table~\ref{tab:capacity} summarizes four main
differences between the cases.

\subsubsection{PSH displaces BESS in the optimized portfolio}
When both technologies compete, the optimal plan builds 31.2~GW of PSH and zero stationary batteries. This occurs even though BESS has a lower
per-kW cost and up to 40~GW was available under the same deployment caps as
in the no-PSH case. The energy dimension explains this result. Per kW, the
4-h BESS carries only 2.4~h of usable storage after applying the 20\%
state-of-charge reserve at both ends, one fifth of the usable energy in the
12-h PSH class, which costs 41\% more per kW-year. Portfolio-wide, the
selected PSH fleet holds $\approx$755~GWh of usable storage, compared with
91~GWh usable (152~GWh nominal) for the 38-GW BESS fleet in the no-PSH
plan. The fleets have comparable power ratings but differ eightfold in
usable energy.

\subsubsection{Selected PSH durations}
\begin{figure}[!t]
\centering
\includegraphics[width=\figw]{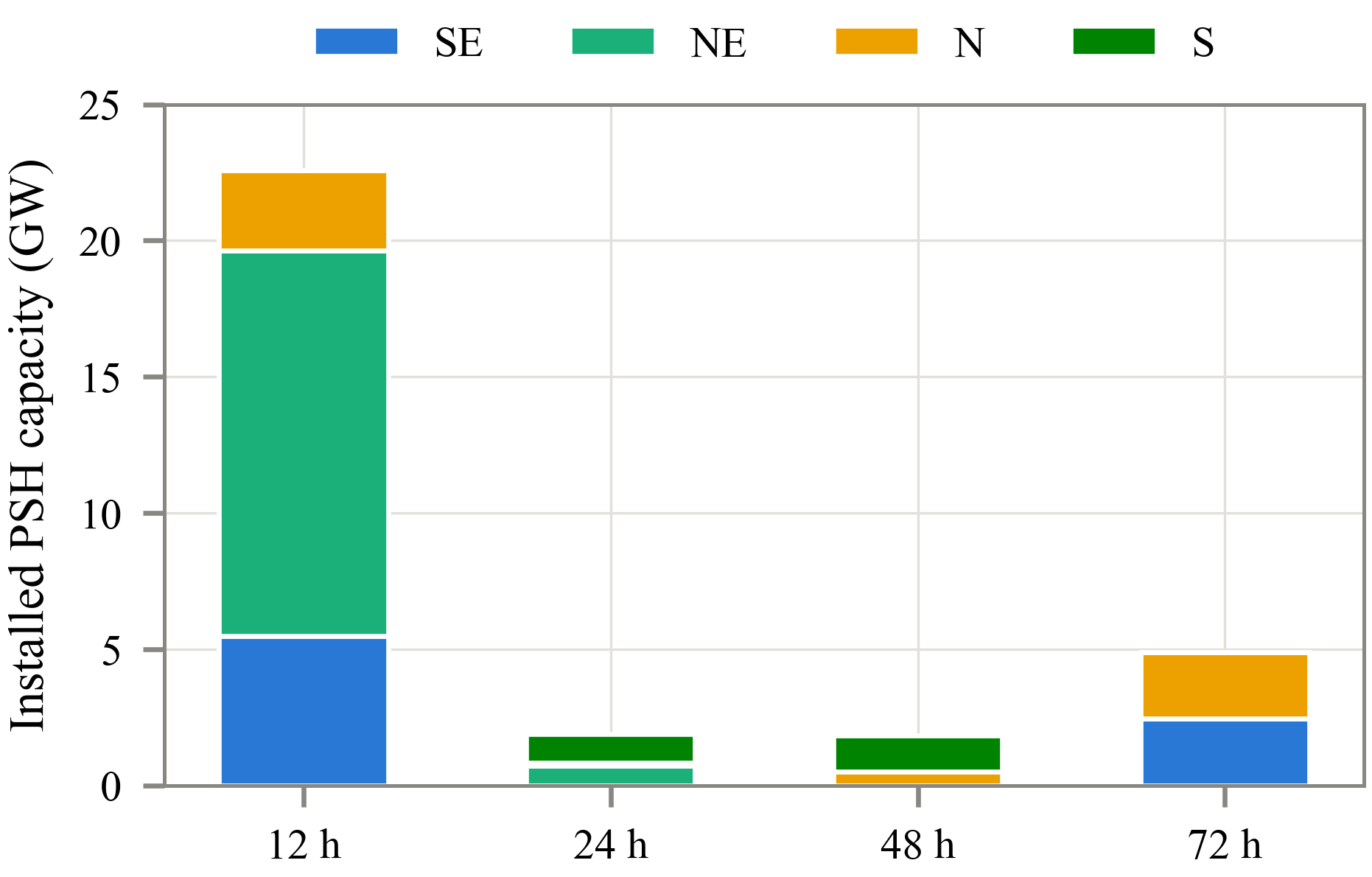}
\caption{Optimal 2050 PSH portfolio by storage duration and subsystem
(with-PSH case; total 31.2\,GW / 755\,GWh). Duration classes 4--8\,h and
100--144\,h receive no investment.}
\label{fig:pshdur}
\end{figure}

Fig.~\ref{fig:pshdur} shows the 2050 PSH portfolio by duration class and
subsystem. The 12-h class dominates (22.6~GW, led by 14.1~GW in the
Northeast, where it firms the wind hub), but the plan also selects 1.9~GW of
24-h, 1.8~GW of 48-h, and 4.9~GW of 72-h units. These longer-duration
classes bridge multi-day wind lulls and the seasonal must-run windows of the
LNG-linked thermal fleet. Neither the
4--8-h nor the 100--144-h classes are selected: the former are dominated by
the 12-h class (only 10\% costlier than 4-h PSH), while the latter's cost
premium is not yet justified at 2050 demand levels.

\subsubsection{Changes in thermal and solar capacity}
Relative to the no-PSH plan, PSH avoids 7.8~GW of new gas-fired capacity and
2.7~GW of thermal re-contracting, while \emph{increasing} optimal solar PV
build-out by 15.6~GW (61.5 vs.\ 45.9~GW in 2050). Low-cost midday solar
energy becomes more valuable when it can be shifted into the evening and
night.

The avoided gas capacity is not distributed evenly across the three candidate
families. Open-cycle peakers account for 5.4~GW of it: the with-PSH plan
builds 2.4~GW of peaking capacity, all of it in the Southeast, against
7.8~GW spread over all four subsystems without PSH. LNG combined cycle
accounts for the remaining 2.4~GW (16.6 vs.\ 19.0~GW). The 10~GW pre-salt
combined-cycle block is built in both plans and in the same year, because its
must-run obligation and low variable cost make it economic whatever the
storage set. Storage therefore substitutes for the capacity that would
otherwise cover short peaks, and for part of the mid-merit LNG fleet, but not
for inflexible must-run generation.

Onshore wind (122~GW) and interconnection reinforcements
($\approx$19~GW, dominated by the 17.3-GW SE--NE corridor expansion) are
nearly identical across cases, indicating that
PSH acts primarily as a substitute for thermal capacity and batteries rather
than for transmission.

\subsubsection{Timing of PSH and corridor investment}
\begin{figure}[!t]
\centering
\includegraphics[width=\figw]{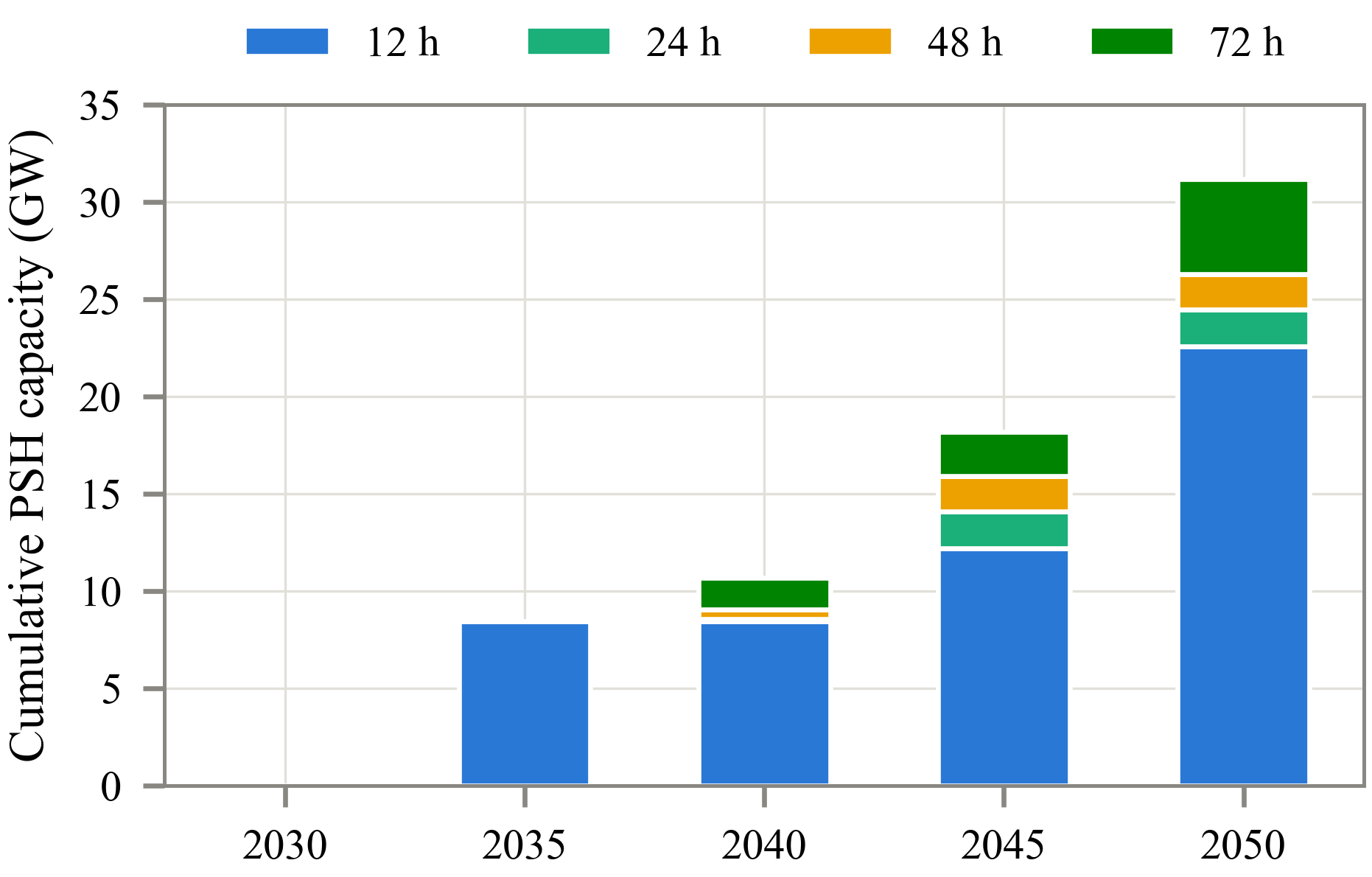}
\caption{PSH build trajectory by duration class (with-PSH case, cumulative
GW). Long-duration classes (48--72\,h) enter in 2040; the 2050 total is
31.2\,GW.}
\label{fig:traj}
\end{figure}

Fig.~\ref{fig:traj} shows the PSH build trajectory by duration class, and the
thermal rows of Table~\ref{tab:capacity} identify the capacity that each
stage replaces. PSH enters at scale in 2035, with 8.4~GW of 12-h capacity in
the Southeast and North, the load centers where it firms the first large
solar tranche. In that year the with-PSH plan builds no open-cycle peakers,
against 3.1~GW in the no-PSH plan, so the earliest PSH blocks substitute for
peaking capacity. The long-duration classes (48--72~h) appear in 2040, the
year in which the two plans first diverge in LNG combined cycle
(1.0 vs.\ 2.6~GW) and inflexible LNG-linked capacity starts setting the
seasonal pattern. The single largest addition is the 2045--2050 block of 12-h
PSH in the Northeast (14.1~GW by 2050), which follows the completion of the
17.3-GW SE--NE corridor expansion in 2045: with the corridor in place,
Northeast storage can arbitrage between the wind hub and the Southeast load
center, and the model concentrates storage where most renewable energy is
produced. Each planning year is optimized without foresight of later years,
so the trajectory records the year in which each addition becomes economic
rather than a coordinated build sequence.

\subsection{System costs}
\label{sec:results:costs}

\begin{table}[!t]
\caption{Annualized system cost by planning year (billion US\$/yr).
Components may not sum exactly to totals because of rounding.}
\label{tab:costs}
\centering
\small
\begin{tabular}{@{}lrrr@{\hskip 10pt}rrr@{\hskip 10pt}r@{}}
\toprule
& \multicolumn{3}{c}{With PSH} & \multicolumn{3}{c}{No PSH} & \\
\cmidrule(r{10pt}){2-4}\cmidrule(r{10pt}){5-7}
Year & CAPEX & OPEX & Total & CAPEX & OPEX & Total & $\Delta$Total \\
\midrule
2030 &       &   2.6 &   2.6 &       &   2.6 &   2.6 & $   $ \\
2035 &   6.7 &   5.3 &  12.0 &   6.8 &   6.1 &  12.9 & $-0.9$ \\
2040 &  17.7 &   5.3 &  23.0 &  17.2 &   7.0 &  24.2 & $-1.2$ \\
2045 &  35.3 &  13.7 &  49.0 &  33.4 &  17.7 &  51.2 & $-2.2$ \\
2050 &  50.9 &  17.5 &  68.3 &  50.5 &  22.9 &  73.3 & $-5.0$ \\
\bottomrule
\end{tabular}
\end{table}

Table~\ref{tab:costs} reports the annualized cost trajectory (cumulative
investment annuities plus expected yearly operating cost). Availability of
PSH reduces the total 2050 system cost by US\$~5.0~billion per year
($-6.8\%$). Annualized CAPEX is 0.7\% higher. The
US\$~6.1~billion/yr PSH annuity is almost offset by avoided BESS, gas, and
thermal re-contracting annuities and by a less costly, solar-heavier
renewable mix. Operating cost falls by US\$~5.4~billion/yr ($-23.5\%$).
In 2050, avoided curtailment (9.5~TWh at 30~US\$/MWh) and avoided unserved
energy (358~GWh at the VOLL) account for US\$~0.3 and
US\$~0.6~billion/yr, respectively. These two terms represent 16\% of the
operating-cost difference; the remaining US\$~4.5~billion/yr is avoided
thermal variable cost. PSH value in the early years therefore comes mainly from
absorbing renewable surplus rather than displacing fuel. Across the horizon
the annual benefit grows from US\$~0.9~billion in 2035 to US\$~5.0~billion in
2050, or 6.8\%, 5.1\%, 4.3\%, and 6.8\% of the no-PSH total in the four
planning years.

\begin{figure}[!t]
\centering
\includegraphics[width=\figw]{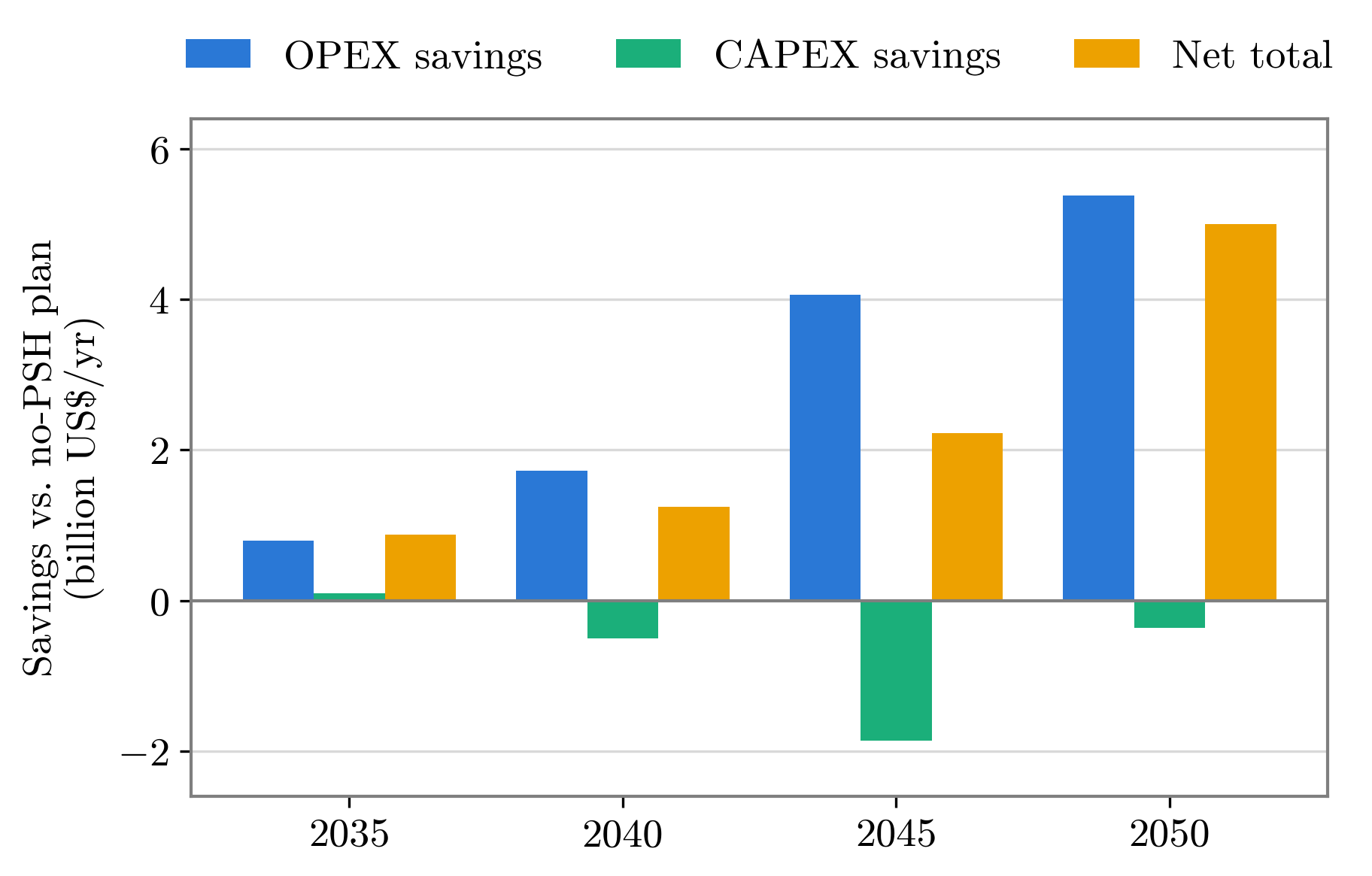}
\caption{Decomposition of the annualized savings of the with-PSH plan
relative to the no-PSH plan. Negative CAPEX savings mean the with-PSH plan
invests more in that year; operating savings more than offset the
difference.}
\label{fig:savings}
\end{figure}

Fig.~\ref{fig:savings} decomposes the savings. OPEX savings dominate in
every year and increase over the horizon. From 2040 onward, the with-PSH
plan spends more on investment annuities, with a maximum difference of
US\$~1.9~billion/yr in 2045. Operating savings exceed this additional
investment by factors of 3.4 in 2040 and 2.2 in 2045, and by more than an
order of magnitude in 2050. The plan thus exchanges higher capital
expenditure for lower operating expenditure.

\subsection{Operational value}
\label{sec:results:ops}

\begin{table*}[!t]
\caption{Expected operational metrics (scenario average).}
\label{tab:ops}
\centering
\small

\begin{tabular}{@{}llrrrr@{}}
\toprule
Metric & Case & 2035 & 2040 & 2045 & 2050 \\
\midrule
Curtailment (TWh)        & PSH    & 14.5 & 12.9 & 16.2 & 20.4 \\
                         & no PSH & 38.2 & 32.1 & 32.5 & 29.9 \\
Curtailment (\% avail.)  & PSH    & 3.2  & 2.2  & 2.2  & 2.3 \\
                         & no PSH & 8.4  & 5.5  & 4.6  & 3.5 \\
Thermal generation (TWh)& PSH    & 80.3 & 98.4 & 216.4 & 263.8 \\
                         & no PSH & 83.6 & 110.1 & 248.0 & 297.9 \\
Energy not served (GWh) & PSH    & 34.5 & 137.6 & 173.3 & 148.4 \\
                         & no PSH & 46.5 & 342.6 & 528.5 & 506.0 \\
LRMC (US\$/MWh)\textsuperscript{a}
                         & PSH    &  66.9 & 118.2 & 117.7 & 101.2 \\
                         & no PSH &  76.5 & 124.9 & 131.3 & 126.5 \\
\bottomrule
\end{tabular}

\vspace{2pt}

\begin{minipage}{0.75\textwidth}
\footnotesize
\textsuperscript{a} Demand-weighted average over subsystems, hours, and
scenarios. Because the model jointly optimizes investment and operational
decisions, the dual variable associated with the load-balance constraint
represents the long-run marginal cost (LRMC) of meeting an incremental
unit of demand.
\end{minipage}

\end{table*}

\begin{figure}[!t]
\centering
\includegraphics[width=\figw]{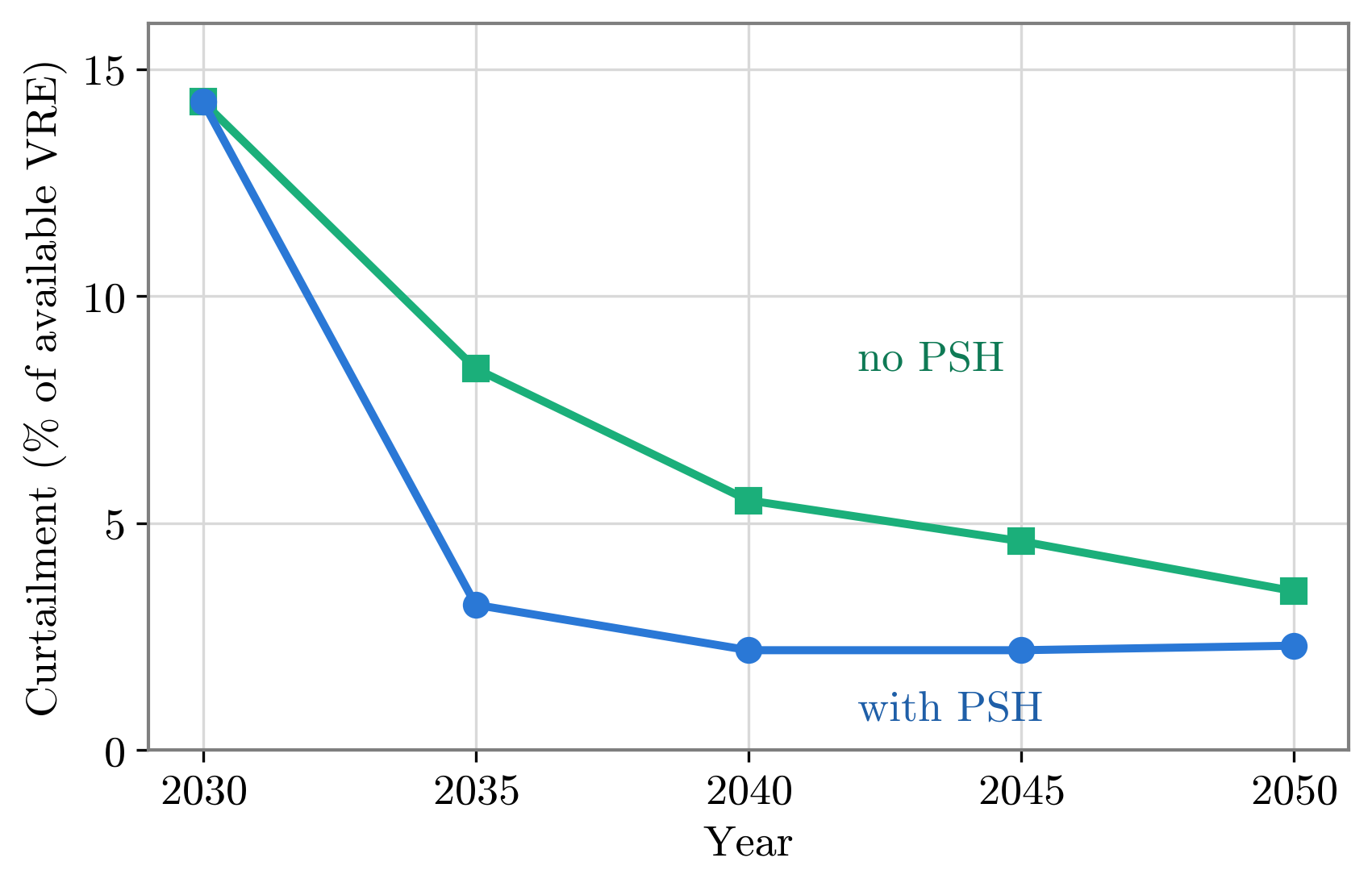}
\caption{Renewable curtailment as a share of available VRE (scenario
average). The 2030 point precedes any storage investment and is common to
both cases.}
\label{fig:curt}
\end{figure}

\begin{figure}[!t]
\centering
\includegraphics[width=\figw]{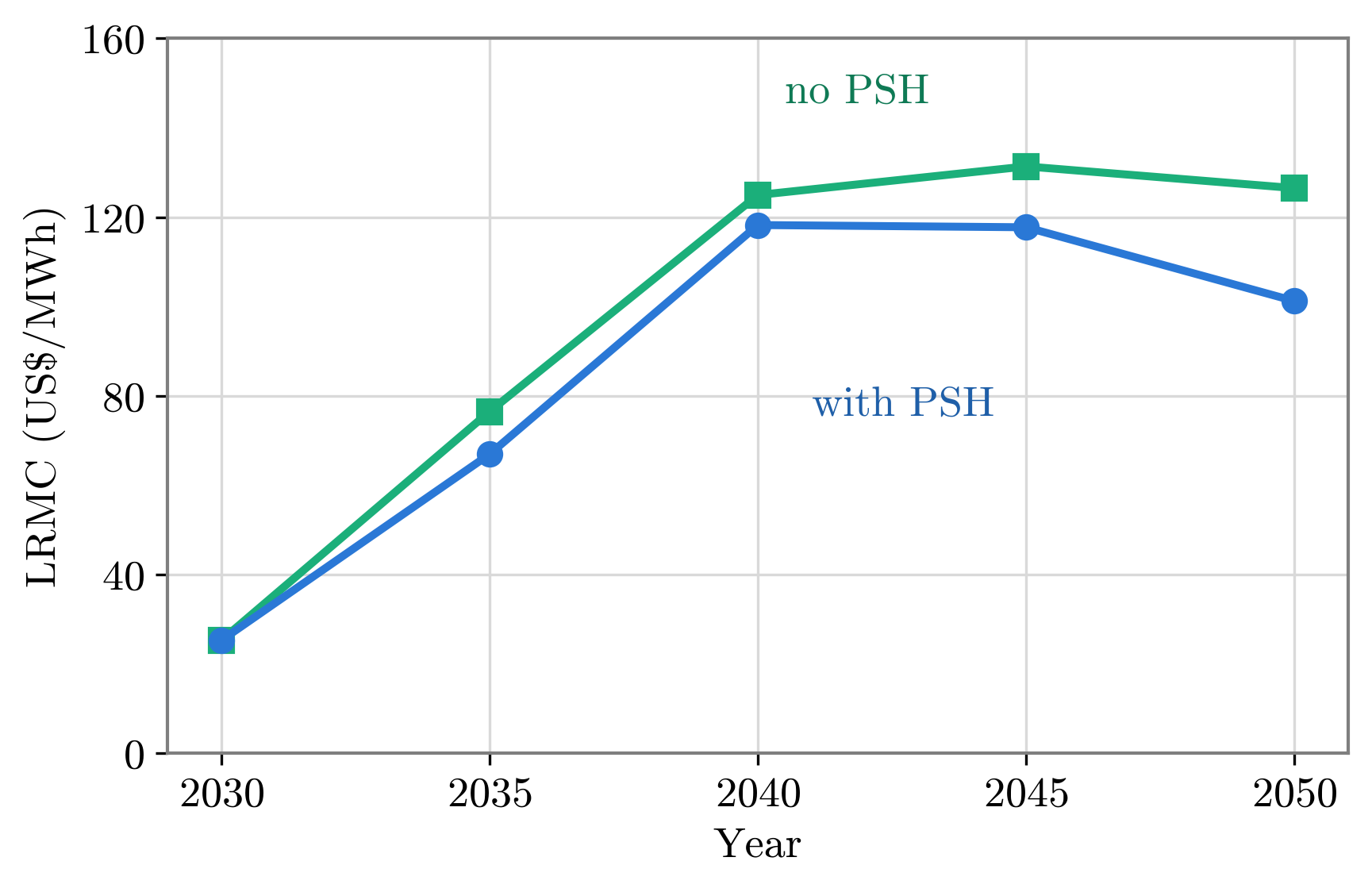}
\caption{Demand-weighted long-run marginal cost (LRMC) of energy.}
\label{fig:cmg}
\end{figure}

Table~\ref{tab:ops} and Figs.~\ref{fig:curt}--\ref{fig:cmg} summarize
operations. In 2035, the no-PSH system curtails 8.4\% of available VRE
(38.2~TWh), compared with 3.2\% (14.5~TWh) with PSH, a 62\% reduction. The difference narrows by 2050 as the no-PSH plan deploys its large BESS fleet.
With PSH, curtailment remains between 2.2\% and 3.2\% throughout the
horizon. Thermal generation falls by 34.1~TWh in 2050 ($-11.4\%$), and
expected energy not served falls by 71\%. The no-PSH value is only 0.03\%
of demand, so this metric indicates episodic adequacy stress rather than
widespread shortage. The demand-weighted LRMC of energy in 2050
falls from 126 to 101~US\$/MWh ($-20\%$).

\begin{table}[!t]
\caption{Mean LRMC of energy by subsystem (US\$/MWh).}
\label{tab:cmgsys}
\centering
\small
\begin{tabular}{@{}llrrrr@{}}
\toprule
Subsystem & Case & 2035 & 2040 & 2045 & 2050 \\
\midrule
Southeast/CW & PSH    &  66.4 & 116.4 & 115.4 &  99.2 \\
             & no PSH &  74.8 & 120.9 & 126.4 & 121.7 \\
South        & PSH    &  67.1 & 117.8 & 116.5 &  98.9 \\
             & no PSH &  75.5 & 122.2 & 127.7 & 118.5 \\
Northeast    & PSH    &  60.3 & 105.6 & 108.1 &  92.3 \\
             & no PSH &  68.8 & 112.4 & 120.4 & 116.9 \\
North        & PSH    &  63.0 & 116.0 & 116.0 &  99.6 \\
             & no PSH &  71.6 & 120.1 & 126.7 & 122.0 \\
\bottomrule
\end{tabular}
\end{table}

Table~\ref{tab:cmgsys} reports the LRMC by subsystem. Every subsystem's LRMC falls by 16--21\% in 2050. The Northeast, which contains the wind hub and the largest PSH block, has the lowest absolute prices in both cases and the largest PSH-related reduction (21\%
in 2050, compared with 16--19\% elsewhere). The SE--NE spread is therefore
\emph{wider} with PSH in every planning year (6.9 vs.\ 4.8~US\$/MWh in
2050). The thermal rows of Table~\ref{tab:capacity} account for the
asymmetry. Most of the 12-h fleet is built in the Northeast (14.1 of
22.6~GW), alongside the wind hub, and it removes that subsystem's peaking
requirement entirely: by 2050 the no-PSH plan has built 3.1~GW of open-cycle
capacity in the Northeast, the with-PSH plan none. The Southeast is the only
subsystem that still builds peakers when PSH is available (2.4~GW, against
2.8~GW without it), and it carries 55\% of national demand behind an import
corridor expanded by nearly the same amount in both cases, so thermal
generation continues to set its evening and dry-season prices even as the
additional 15.6~GW of solar lowers its midday cost. The corridor does not
arbitrage the resulting difference away. Price and capacity results therefore
point the same way: PSH substitutes for thermal
capacity and batteries rather than for transmission.

\begin{figure}[!t]
\centering
\includegraphics[width=\figw]{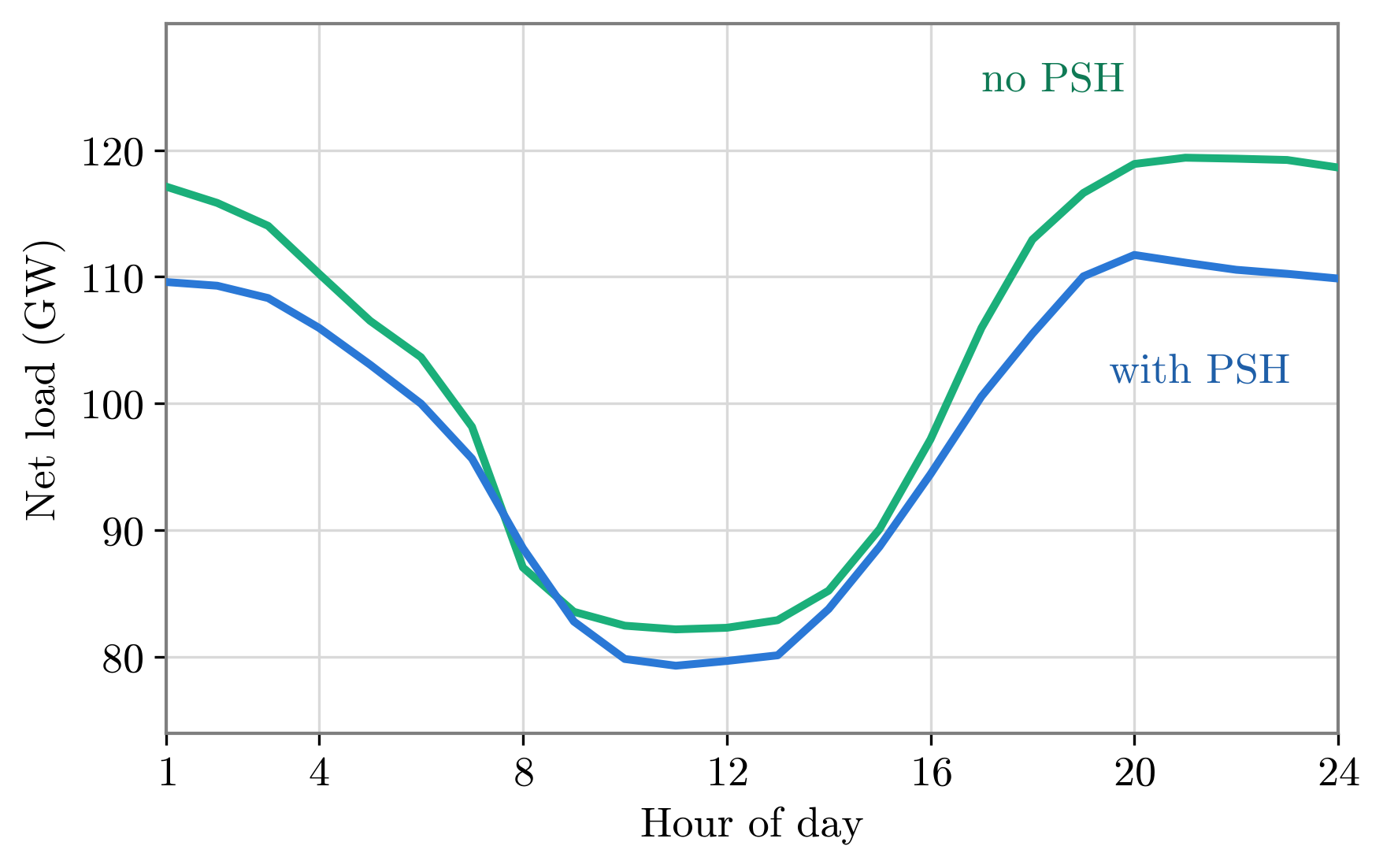}
\caption{Average daily net load after storage action in 2050 (demand minus
renewable generation minus net storage injection; mean over all days and
scenarios). The 31.2-GW PSH fleet holds the residual net load within a
32.4-GW band (79.3--111.7\,GW), versus a 37.2-GW band for the 38.0-GW BESS
fleet
of the storage-restricted plan, even though the with-PSH plan carries
15.6\,GW more solar and has a deeper pre-storage valley.}
\label{fig:daily}
\end{figure}

Storage operation explains these differences. Before storage dispatch, the
with-PSH plan has a 78-GW daily swing in \emph{net load} (demand minus
renewable generation), from a 52-GW midday valley to a 130-GW evening peak.
This exceeds the no-PSH plan's 66-GW swing because the with-PSH plan carries
15.6~GW more solar. Fig.~\ref{fig:daily} shows the average net load after
storage dispatch. The PSH fleet absorbs up to 27~GW of midday surplus and
injects up to 18~GW during the evening peak. Positive injection continues
overnight, and the residual net load seen by the hydro-thermal system spans
32.4~GW (79.3--111.7~GW). The
4-h BESS fleet charges at roughly two thirds of that rate, and its overnight
delivery falls below 3~GW in the early-morning hours. Its residual net load
spans 37~GW despite a shallower initial valley.

\begin{figure}[!t]
\centering
\includegraphics[width=\figw]{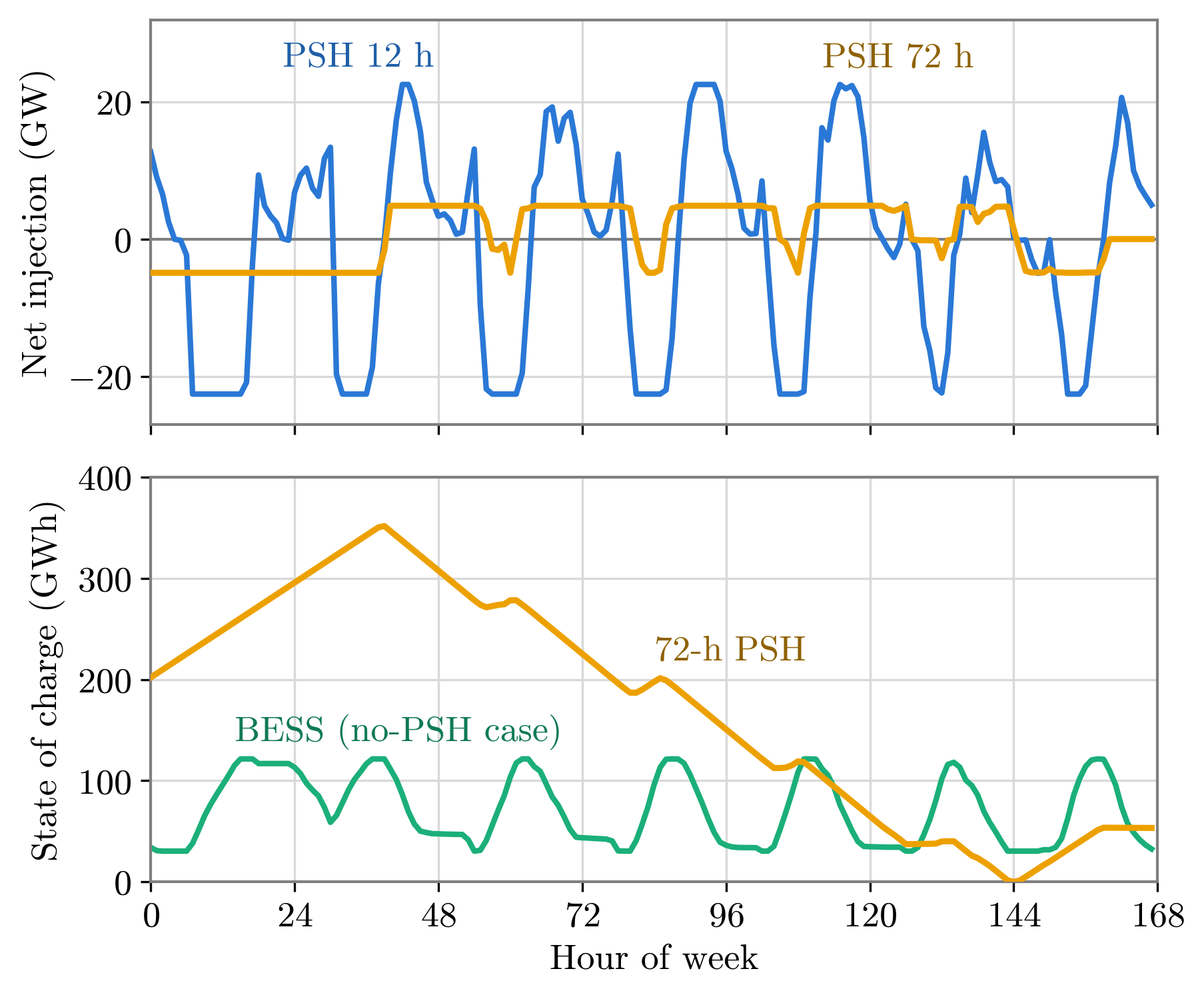}
\caption{Storage operation in an illustrative week of 2050 (scenario~1;
hours 2689--2856). Top: net injection of the 12-h and 72-h PSH classes. Bottom: state of charge of the 72-h PSH fleet (with-PSH
case) and of the entire BESS fleet (no-PSH case). The 72-h units execute a
single full-volume weekly cycle that a 4-h battery portfolio cannot replicate.}
\label{fig:week}
\end{figure}

Weekly operation distinguishes the duration classes.
Fig.~\ref{fig:week} plots an illustrative April week of 2050 (scenario~1).
The 12-h class cycles daily, mirroring the solar cycle, while the 72-h fleet
performs one full charge--discharge excursion over the week
($\approx$350~GWh, its entire usable volume). It charges before a low-wind
period and then discharges for several consecutive days. The BESS fleet in
the no-PSH case is bounded by its 91-GWh usable volume and cannot provide
the same multi-day service. Gas thermal units meet more of the residual
balancing requirement, contributing to the OPEX difference in
Table~\ref{tab:costs}. Total 2050 storage throughput is similar in the two
plans (71.4~TWh discharged by PSH vs.\ 75.9~TWh by BESS), and the realized
round-trip efficiencies differ by only three percentage points (80\% vs.\
83\%). Usable energy per cycle differs much more. Measured against usable
energy capacity, the BESS fleet performs $\approx$830 equivalent full cycles
per year ($\approx$2.3 per day), compared with $\approx$95 for PSH. Each
PSH cycle moves roughly eight times more energy.

\subsection{Seasonal patterns}
\label{sec:results:seasonal}

\begin{figure}[!t]
\centering
\includegraphics[width=\figw]{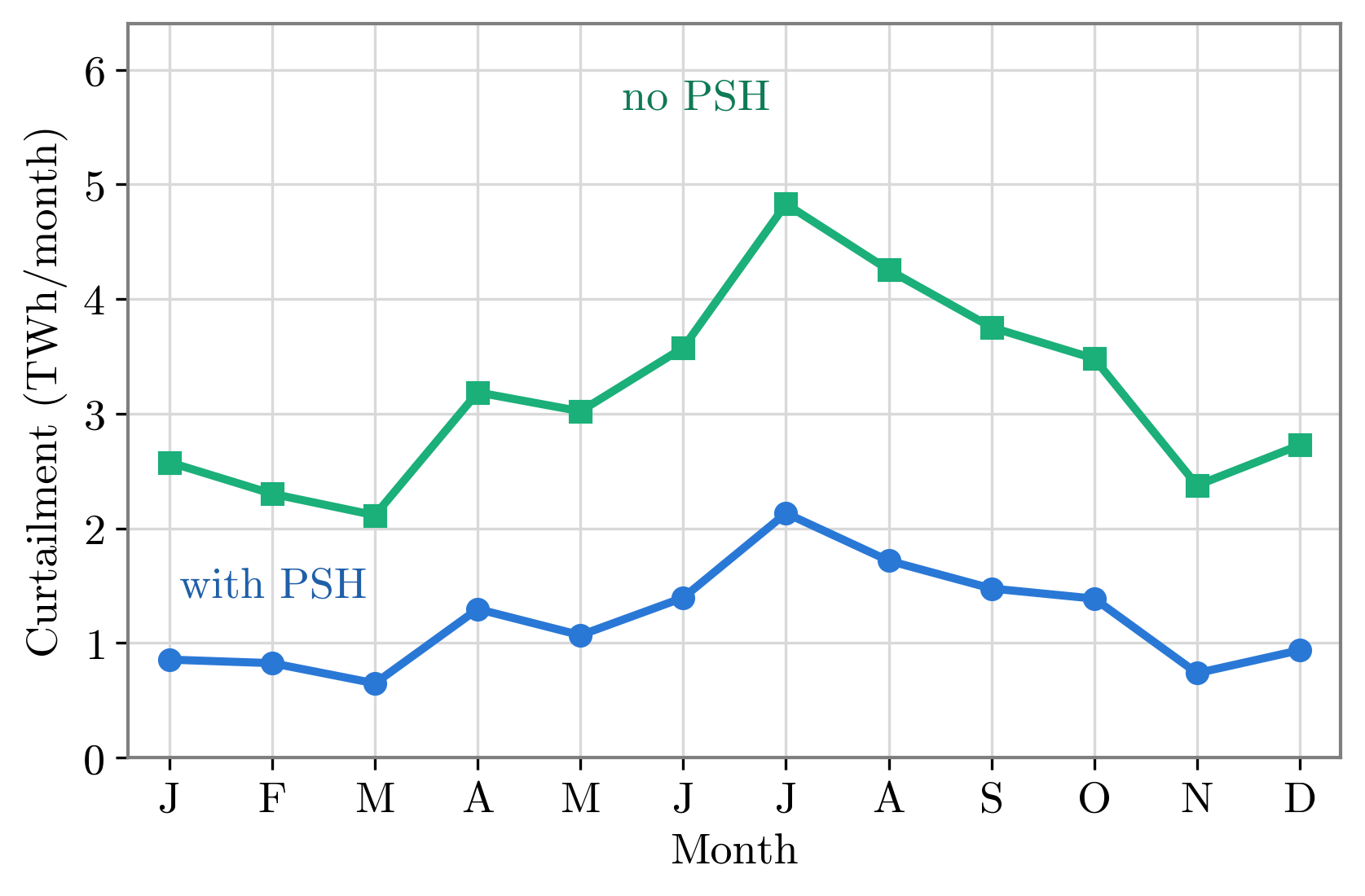}
\caption{Monthly renewable curtailment in 2035 (scenario average). The gap
peaks in the June--October dry season, when Northeast wind output is
highest and LNG-linked thermal plants are must-run.}
\label{fig:monthly}
\end{figure}

Monthly results show how PSH interacts with the thermal contract regimes.
Fig.~\ref{fig:monthly} reports curtailment in 2035, the year with the largest
relative difference. Without PSH, the system curtails 2.1--4.8~TWh each
month, with a July peak during the June--November period when LNG-contracted
combined-cycle plants are must-run and Northeast wind output is strongest.
With PSH, monthly curtailment remains below 2.2~TWh. In 2050, PSH discharge
varies from 5.1 to 6.8~TWh/month, and monthly thermal generation is
2.2--3.2~TWh lower than in the no-PSH case in every month. In the modeled
scenarios, savings arise throughout the year through daily solar
shifting in the wet season and multi-day wind and must-run balancing in the
dry season.

\subsection{Computational performance}
\label{sec:results:comp}

Each planning year is a deterministic-equivalent LP with all ten scenarios
at hourly resolution; the largest instances reach $5.2\times10^{7}$
constraints, $4.5\times10^{7}$ variables and $1.6\times10^{8}$ nonzeros,
solved with the Xpress barrier algorithm without crossover on 16 threads
($\le$128~GB RAM) on AWS r7i instance. Adding the 36 PSH candidates
triples total solution time (13.3~h vs.\ 4.4~h for the full five-year
horizon). Studies with larger scenario
sets require scenario-decomposition methods.

\section{Discussion}
\label{sec:discussion}

\subsection{Storage duration and portfolio substitution}

The study is set up as a portfolio question: which storage mix an hourly
expansion model selects when the candidate set spans 4--144-h durations.
Widening the candidate set in this way changes the outcome. The selected
portfolio and the 4-h-only benchmark have comparable power ratings (31.2 and
38.0~GW) but different energy volumes. The PSH fleet provides
approximately 755~GWh of
usable storage, compared with 91~GWh usable (152~GWh nominal) in the
BESS-only case. Both fleets discharge a similar amount of energy over the
year, and their realized round-trip efficiencies differ by only three
percentage points. What separates them is the energy moved per cycle and the
timescale of that cycle: 12-h PSH shifts energy daily, whereas the 24--72-h
classes bridge multi-day events. This result is consistent with studies
showing that storage value depends strongly on the energy-to-power
ratio~\cite{sepulveda2021,denholm2021}.

The comparison does not establish that PSH will displace BESS under all cost
and deployment assumptions. The BESS-only fleet reaches 95\% of its national
deployment limit, and only a 4-h battery class is represented. The result
instead shows that evaluating storage only in GW, or against a single fixed
duration, can produce a materially different capacity mix from a model in
which storage energy and duration are selectable.

\subsection{System value beyond storage arbitrage}

Admitting PSH changes investments outside the storage portfolio. Relative to
the BESS-only case, the with-PSH plan builds 15.6~GW more solar PV, 7.8~GW
less gas-fired capacity (5.4~GW of it open-cycle peaking capacity), and
2.7~GW less thermal re-contracting by 2050. The 10~GW pre-salt
combined-cycle block is unchanged, so the substitution acts on peaking and
mid-merit capacity rather than on must-run generation.
Transmission expansion is nearly unchanged. In this model, PSH
complements the SE--NE corridor and substitutes primarily for
batteries and thermal capacity rather than for transmission.

The cost results also separate system value from project revenue. Annualized
CAPEX is 0.7\% higher with PSH in 2050, while operating cost is 23.5\% lower.
Avoided thermal variable cost accounts for most of the difference, with
additional contributions from lower curtailment and unserved energy. These
benefits accrue across the system and may not be recovered through energy
arbitrage alone~\cite{jacobson2023}. Because the model does not simulate
project revenues or market bidding, it identifies a system-optimal
investment but does not establish project bankability.

The operating results show how that value is created. The
12-h fleet narrows the daily residual-net-load range while accommodating
additional solar capacity. The 72-h fleet makes a full weekly energy
excursion during a multi-day wind event, and monthly savings persist under
both wet- and dry-season operating conditions. The value attributed to PSH
combines capacity substitution with operation across several
timescales.

\subsection{Planning and procurement implications}

For expansion studies, storage candidate sets should distinguish power from
energy and span the durations relevant to the modeled system. A single
power-capacity target cannot represent the different daily and multi-day
services observed here. For procurement, the result supports specifying
deliverable energy or duration alongside MW capability. Brazil's capacity
reserve auction for storage currently admits batteries but not
PSH~\cite{mme136}. Under the modeled assumptions, that eligibility choice
changes both the storage portfolio and total system cost.

The model does not determine a preferred remuneration mechanism. It does,
however, identify a gap between modeled system value and energy-market
revenue. China's two-part capacity and energy tariff~\cite{iha2024} and
state-level floor-and-ceiling contracts for long-duration storage in
Australia~\cite{aemo_isp} illustrate approaches that remunerate more than
energy arbitrage. Any Brazilian application would still require
technology-neutral qualification rules and site-specific environmental and
grid assessments for closed-loop PSH.

\subsection{Limitations}

Six modeling choices bound the interpretation of the results.
\begin{enumerate}
  \item The four-subsystem representation omits intra-regional network
  effects. The finding that PSH does not displace transmission therefore
  applies only at the corridor level. Named-site analysis with GIS screening
  and explicit connection costs may identify additional transmission value.

  \item The rolling horizon is myopic: each planning year is optimized
  without foresight of later years, and no inter-year discounting is applied.
  Costs are annual snapshots rather than a present-value total. Capacity also
  materializes without construction lead times, which may favor PSH.

  \item Perfect foresight within each 8{,}760-hour scenario may overstate
  storage value, especially for the multi-day pre-charging shown in
  Fig.~\ref{fig:week}. Operations also omit unit commitment, while
  retirements are exogenous.

  \item The BESS alternative is limited to a 4-h class with a 60\% usable
  state-of-charge window, 83\% round-trip efficiency, and a deployment cap.
  Longer-duration batteries, different build limits, or steeper cost declines
  could weaken the crowding-out result. Sensitivities to fuel prices, carbon
  prices, and the curtailment penalty are not evaluated even though the
  estimated benefit is dominated by operating cost.

  \item Investment variables are continuous and do not represent project
  lumpiness.

  \item The linear storage formulation does not explicitly prohibit
  simultaneous charging and discharging. Efficiency losses and positive
  throughput costs discourage this behavior but do not provide a mathematical
  guarantee. The model also excludes ancillary-service and inertia value,
  services that variable-speed PSH can provide
  \cite{dong2020,ieee_tr134}.
\end{enumerate}

\section{Conclusion}
\label{sec:conclusion}

This study asks what an expansion model selects when storage candidates of
many durations compete and the duration mix is left to the optimization. The
framework is an hourly,
rolling-horizon, two-stage stochastic
generation--transmission--storage program, and we exercise it on the
Brazilian Interconnected System over 2030--2050. Rather than assigning PSH a
fixed operating role, the model selects its power, energy, duration,
location, and dispatch while it competes with 4-h BESS under otherwise
identical system assumptions.

With PSH available, the model builds 31.2~GW / 755~GWh of PSH and no BESS by
2050. The selected mix contains a dominant 12-h block and 8.6~GW of
24--72-h capacity. The BESS-only case instead builds 38.0~GW / 152~GWh of
batteries, near their deployment limit, and 7.8~GW more gas-fired capacity,
of which 5.4~GW is open-cycle peaking capacity.
The two storage portfolios have comparable power but differ approximately
eightfold in usable energy. Daily cycling of the 12-h fleet and multi-day
cycling of the longer-duration classes explain this difference.

The cost difference is created in operation, not in investment. Annualized
CAPEX differs by less than 1\% between the two plans, while operating cost is
23.5\% lower with PSH, for a total 2050 saving of US\$~5.0~billion/year.
Most of that saving is avoided thermal fuel: the with-PSH plan generates
34~TWh less from thermal units, carries 15.6~GW more solar PV, and settles at
a demand-weighted LRMC 20\% below the BESS-only plan. Curtailment is
60--62\% lower in 2035--2040. Transmission expansion is nearly unchanged.

These results are conditional on the assumed storage costs, perfect foresight within each scenario, and the zonal network representation. Within those bounds, the main finding is that storage
duration should be resolved inside expansion planning rather than fixed
before it. Representing
storage only by MW capacity or by a single 4-h class can materially change
the estimated storage mix and system value. That conclusion is a property of
the modeling choice rather than of the Brazilian system: the quantities
reported here depend on Brazilian costs, resources, and contract structures,
but any system whose residual load has structure beyond the daily cycle can
be expected to select a duration mix that a single-duration candidate set
cannot represent. Further work should test
named-site PSH candidates, longer-duration BESS alternatives, storage-cost
and deployment-limit sensitivities, imperfect foresight, and explicit
ancillary-service requirements.

\section*{CRediT authorship contribution statement}

\textbf{Rafael Benchimol Klausner:} Conceptualization, Methodology,
Software, Formal analysis, Investigation, Visualization, Writing -- original draft. \textbf{Rafael Kelman:} Conceptualization, Methodology, Supervision, Funding acquisition, Writing -- review \& editing.

\section*{Declaration of competing interest}

The authors declare that they have no known competing financial interests or
personal relationships that could have appeared to influence the work
reported in this paper.

\section*{Funding}

This work was funded by the Research and Development (PD\&I) Programme of the
Brazilian Electricity Regulatory Agency (ANEEL) under project
PD-00387-0125, sponsored by CTG Brasil. The funding source had no role
in the study design, in the analysis and interpretation of the results, or in
the decision to submit the article for publication.

\section*{Acknowledgments}

This work was developed within project PD-00387-0125, executed by the
PSR--HPPA--INTT consortium. The authors thank the project team (PSR, HPPA,
INTT) and CTG Brasil for the project framework and data assumptions used in
this study.

\section*{Declaration of generative AI and AI-assisted technologies in the
manuscript preparation process}

During the preparation of this work, the authors used Anthropic Claude to improve
language, readability, and organization. After using this tool, the authors
reviewed and edited the content as needed and take full responsibility for
the content of the published article.

\bibliographystyle{elsarticle-num}
\bibliography{main}

\begin{thebibliography}{10}
\expandafter\ifx\csname url\endcsname\relax
  \def\url#1{\texttt{#1}}\fi
\expandafter\ifx\csname urlprefix\endcsname\relax\def\urlprefix{URL }\fi
\expandafter\ifx\csname href\endcsname\relax
  \def\href#1#2{#2} \def\path#1{#1}\fi

\bibitem{bistline2020}
J.~Bistline, W.~Cole, G.~Damato, J.~DeCarolis, W.~Frazier, V.~Linga, C.~Marcy,
  C.~Namovicz, K.~Podkaminer, R.~Sims, M.~Sukunta, D.~Young, Energy storage in
  long-term system models: a review of considerations, best practices, and
  research needs, Progress in Energy 2~(3) (2020) 032001.
\newblock \href {https://doi.org/10.1088/2516-1083/ab9894}
  {\path{doi:10.1088/2516-1083/ab9894}}.

\bibitem{iha2024}
{International Hydropower Association},
  \href{https://www.hydropower.org/publications/2024-world-hydropower-outlook}{2024
  world hydropower outlook: Opportunities to advance net-zero}, Tech. rep.,
  IHA, London, U.K., accessed 28 July 2026 (2024).
\newline\urlprefix\url{https://www.hydropower.org/publications/2024-world-hydropower-outlook}

\bibitem{blakers2021}
A.~Blakers, M.~Stocks, B.~Lu, C.~Cheng, A review of pumped hydro energy
  storage, Progress in Energy 3~(2) (2021) 022003.
\newblock \href {https://doi.org/10.1088/2516-1083/abeb5b}
  {\path{doi:10.1088/2516-1083/abeb5b}}.

\bibitem{koritarov2022}
V.~Koritarov, et~al.,
  \href{https://publications.anl.gov/anlpubs/2022/05/175341.pdf}{A review of
  technology innovations for pumped storage hydropower}, Tech. Rep. ANL-22/08,
  Argonne National Laboratory, accessed 28 July 2026 (2022).
\newline\urlprefix\url{https://publications.anl.gov/anlpubs/2022/05/175341.pdf}

\bibitem{zhao2024}
Z.~Zhao, et~al., Beyond fixed-speed pumped storage: A comprehensive evaluation
  of different flexible pumped storage technologies in energy systems, Journal
  of Cleaner Production 434 (2024) 139994.
\newblock \href {https://doi.org/10.1016/j.jclepro.2023.139994}
  {\path{doi:10.1016/j.jclepro.2023.139994}}.

\bibitem{papadakis2023}
N.~C. Papadakis, M.~Fafalakis, D.~Katsaprakakis, A review of pumped hydro
  storage systems, Energies 16~(11) (2023) 4516.
\newblock \href {https://doi.org/10.3390/en16114516}
  {\path{doi:10.3390/en16114516}}.

\bibitem{stocks2021}
M.~Stocks, R.~Stocks, B.~Lu, C.~Cheng, A.~Blakers, Global atlas of closed-loop
  pumped hydro energy storage, Joule 5~(1) (2021) 270--284.
\newblock \href {https://doi.org/10.1016/j.joule.2020.11.015}
  {\path{doi:10.1016/j.joule.2020.11.015}}.

\bibitem{krishnan2016}
V.~Krishnan, J.~Ho, B.~F. Hobbs, A.~L. Liu, J.~D. McCalley, M.~Shahidehpour,
  Q.~P. Zheng, Co-optimization of electricity transmission and generation
  resources for planning and policy analysis: review of concepts and modeling
  approaches, Energy Systems 7~(2) (2016) 297--332.
\newblock \href {https://doi.org/10.1007/s12667-015-0158-4}
  {\path{doi:10.1007/s12667-015-0158-4}}.

\bibitem{jacobson2023}
M.~Jacobson, J.~Tan, E.~Muljadi, D.~Corbus, Z.~Dong, J.~Kim, E.~Bailey,
  M.~Pevarnik, M.~Racine, A.~{St.\ Hilaire}, C.~Hodge,
  \href{https://docs.nrel.gov/docs/fy23osti/74775.pdf}{An assessment of
  deploying advanced pumped storage hydropower technology in {U.S.} electricity
  markets}, Tech. Rep. NREL/TP-7A40-74775, National Renewable Energy
  Laboratory, Golden, CO, accessed 28 July 2026 (Jul. 2023).
\newline\urlprefix\url{https://docs.nrel.gov/docs/fy23osti/74775.pdf}

\bibitem{jump}
M.~Lubin, O.~Dowson, J.~D. Garcia, J.~Huchette, B.~Legat, J.~P. Vielma, {JuMP}
  1.0: Recent improvements to a modeling language for mathematical
  optimization, Mathematical Programming Computation 15 (2023) 581--589.
\newblock \href {https://doi.org/10.1007/s12532-023-00239-3}
  {\path{doi:10.1007/s12532-023-00239-3}}.

\bibitem{xpress}
{FICO}, \href{https://www.fico.com/en/products/fico-xpress-optimization}{{FICO
  Xpress Optimizer}, version 9.8}, Fair Isaac Corporation, accessed 28 July
  2026 (2026).
\newline\urlprefix\url{https://www.fico.com/en/products/fico-xpress-optimization}

\bibitem{iea_brazil2025}
{International Energy Agency},
  \href{https://www.iea.org/reports/brazil-2025}{Energy policy review: {Brazil}
  2025}, Tech. rep., IEA, Paris, France, accessed 30 July 2026 (2025).
\newline\urlprefix\url{https://www.iea.org/reports/brazil-2025}

\bibitem{epe_pde2034}
{Empresa de Pesquisa Energ\'etica},
  \href{https://www.epe.gov.br/pt/publicacoes-dados-abertos/publicacoes/plano-decenal-de-expansao-de-energia-2034}{Plano
  decenal de expans\~ao de energia 2034}, Tech. rep., EPE / Minist\'erio de
  Minas e Energia, Rio de Janeiro, Brazil, accessed 28 July 2026 (2025).
\newline\urlprefix\url{https://www.epe.gov.br/pt/publicacoes-dados-abertos/publicacoes/plano-decenal-de-expansao-de-energia-2034}

\bibitem{mme136}
{MME},
  \href{https://www.epe.gov.br/pt/imprensa/noticias/lrcap-armazenamento-2026-epe-publica-orientacoes-para-cadastramento-dos-projetos}{Portaria
  normativa n\textsuperscript{o} 136/gm/mme, de 1\textsuperscript{o} de junho
  de 2026: diretrizes dos leil\~oes de reserva de capacidade armazenamento de
  2026}, Minist\'erio de Minas e Energia, accessed 28 July 2026 (2026).
\newline\urlprefix\url{https://www.epe.gov.br/pt/imprensa/noticias/lrcap-armazenamento-2026-epe-publica-orientacoes-para-cadastramento-dos-projetos}

\bibitem{poncelet2016}
K.~Poncelet, E.~Delarue, D.~Six, J.~Duerinck, W.~D'haeseleer, Impact of the
  level of temporal and operational detail in energy-system planning models,
  Applied Energy 162 (2016) 631--643.
\newblock \href {https://doi.org/10.1016/j.apenergy.2015.10.100}
  {\path{doi:10.1016/j.apenergy.2015.10.100}}.

\bibitem{helisto2019}
N.~Helist\"o, J.~Kiviluoma, H.~Holttinen, J.~D. Lara, B.-M. Hodge, Including
  operational aspects in the planning of power systems with large amounts of
  variable generation: A review of modeling approaches, WIREs Energy and
  Environment 8~(5) (2019) e341.
\newblock \href {https://doi.org/10.1002/wene.341}
  {\path{doi:10.1002/wene.341}}.

\bibitem{pypsa2018}
T.~Brown, J.~H\"orsch, D.~Schlachtberger, {PyPSA}: {Python} for power system
  analysis, Journal of Open Research Software 6~(1) (2018) 4.
\newblock \href {https://doi.org/10.5334/jors.188}
  {\path{doi:10.5334/jors.188}}.

\bibitem{switch2019}
J.~Johnston, R.~Henriquez-Auba, B.~Maluenda, M.~Fripp, {Switch} 2.0: A modern
  platform for planning high-renewable power systems, SoftwareX 10 (2019)
  100251.
\newblock \href {https://doi.org/10.1016/j.softx.2019.100251}
  {\path{doi:10.1016/j.softx.2019.100251}}.

\bibitem{sepulveda2021}
N.~A. Sepulveda, J.~D. Jenkins, A.~Edington, D.~S. Mallapragada, R.~K. Lester,
  The design space for long-duration energy storage in decarbonized power
  systems, Nature Energy 6 (2021) 506--516.
\newblock \href {https://doi.org/10.1038/s41560-021-00796-8}
  {\path{doi:10.1038/s41560-021-00796-8}}.

\bibitem{denholm2021}
P.~Denholm, W.~Cole, N.~Blair,
  \href{https://docs.nrel.gov/docs/fy23osti/85878.pdf}{Moving beyond 4-hour
  li-ion batteries: Challenges and opportunities for long(er)-duration energy
  storage}, Tech. Rep. NREL/TP-6A40-85878, National Renewable Energy
  Laboratory, Golden, CO, accessed 28 July 2026 (2023).
\newline\urlprefix\url{https://docs.nrel.gov/docs/fy23osti/85878.pdf}

\bibitem{dowling2020}
J.~A. Dowling, K.~Z. Rinaldi, T.~H. Ruggles, S.~J. Davis, M.~Yuan, F.~Tong,
  N.~S. Lewis, K.~Caldeira, Role of long-duration energy storage in variable
  renewable electricity systems, Joule 4~(9) (2020) 1907--1928.
\newblock \href {https://doi.org/10.1016/j.joule.2020.07.007}
  {\path{doi:10.1016/j.joule.2020.07.007}}.

\bibitem{birge2011}
J.~R. Birge, F.~Louveaux, Introduction to Stochastic Programming, 2nd Edition,
  Springer, New York, NY, 2011.
\newblock \href {https://doi.org/10.1007/978-1-4614-0237-4}
  {\path{doi:10.1007/978-1-4614-0237-4}}.

\bibitem{pereira1991}
M.~V.~F. Pereira, L.~M. V.~G. Pinto, Multi-stage stochastic optimization
  applied to energy planning, Mathematical Programming 52 (1991) 359--375.
\newblock \href {https://doi.org/10.1007/BF01582895}
  {\path{doi:10.1007/BF01582895}}.

\bibitem{barbour2016}
E.~Barbour, I.~A.~G. Wilson, J.~Radcliffe, Y.~Ding, Y.~Li, A review of pumped
  hydro energy storage development in significant international electricity
  markets, Renewable and Sustainable Energy Reviews 61 (2016) 421--432.
\newblock \href {https://doi.org/10.1016/j.rser.2016.04.019}
  {\path{doi:10.1016/j.rser.2016.04.019}}.

\bibitem{dong2020}
Z.~Dong, et~al., Developing of quaternary pumped storage hydropower for dynamic
  studies, IEEE Transactions on Sustainable Energy 11~(4) (2020) 2870--2878.
\newblock \href {https://doi.org/10.1109/TSTE.2020.2980585}
  {\path{doi:10.1109/TSTE.2020.2980585}}.

\bibitem{ieee_tr134}
{IEEE PES Power System Dynamic and Performance Committee}, Advanced pumped
  storage hydropower modeling, Tech. Rep. PES-TR134, IEEE Power \& Energy
  Society (Jul. 2025).

\bibitem{hunt2020}
J.~D. Hunt, E.~Byers, Y.~Wada, S.~Parkinson, D.~E. H.~J. Gernaat, S.~Langan,
  D.~P. van Vuuren, K.~Riahi, Global resource potential of seasonal pumped
  hydropower storage for energy and water storage, Nature Communications 11
  (2020) 947.
\newblock \href {https://doi.org/10.1038/s41467-020-14555-y}
  {\path{doi:10.1038/s41467-020-14555-y}}.

\bibitem{cohen2023}
S.~Cohen, V.~Ramasamy, D.~Inman,
  \href{https://docs.nrel.gov/docs/fy23osti/84875.pdf}{A component-level
  bottom-up cost model for pumped storage hydropower}, Tech. Rep.
  NREL/TP-6A40-84875, National Renewable Energy Laboratory, Golden, CO,
  accessed 28 July 2026 (2023).
\newline\urlprefix\url{https://docs.nrel.gov/docs/fy23osti/84875.pdf}

\bibitem{li2023}
J.~Li, Z.~Zhao, D.~Xu, P.~Li, Y.~Liu, M.~A. Mahmud, D.~Chen, The potential
  assessment of pump hydro energy storage to reduce renewable curtailment and
  {CO$_2$} emissions in {Northwest China}, Renewable Energy 212 (2023) 82--96.
\newblock \href {https://doi.org/10.1016/j.renene.2023.04.132}
  {\path{doi:10.1016/j.renene.2023.04.132}}.

\bibitem{zhu2026}
Z.~Zhu, H.~Mao, S.~Zhang, X.~He, D.~Zhang, Spatially resolved modeling of
  pumped storage and hydropower for {China's} carbon neutrality, Energy \&
  Environmental Science 19~(3) (2026) 906--925.
\newblock \href {https://doi.org/10.1039/D5EE05948H}
  {\path{doi:10.1039/D5EE05948H}}.

\bibitem{silva2025}
L.~J. da~Silva, V.~Parente, J.~O.~N. de~Jesus, K.~P.~O. Esquerre, O.~Sahin,
  W.~C. de~Araujo, Pumped hydro storage in the {Brazilian} power industry: A
  sustainable approach to expanding renewable energy, Sustainability 17~(5)
  (2025) 1911.
\newblock \href {https://doi.org/10.3390/su17051911}
  {\path{doi:10.3390/su17051911}}.

\bibitem{aemo_isp}
{Australian Energy Market Operator},
  \href{https://aemo.com.au/energy-systems/major-publications/integrated-system-plan-isp}{Draft
  2026 integrated system plan}, Tech. rep., AEMO, Melbourne, Australia,
  accessed 28 July 2026 (Dec. 2025).
\newline\urlprefix\url{https://aemo.com.au/energy-systems/major-publications/integrated-system-plan-isp}

\bibitem{weber2024}
N.~{de Assis Brasil Weber}, J.~D. Hunt, B.~Zakeri, P.~S. Schneider, F.~S.~A.
  Parente, A.~D. Marques, A.~O. {Pereira Junior}, Seasonal pumped hydropower
  storage role in responding to climate change impacts on the {Brazilian}
  electrical sector, Journal of Energy Storage 87 (2024) 111249.
\newblock \href {https://doi.org/10.1016/j.est.2024.111249}
  {\path{doi:10.1016/j.est.2024.111249}}.

\bibitem{albuquerque2026}
L.~R. Albuquerque, R.~Kelman, T.~Castro, A.~Petrungaro, T.~Andrade,
  \href{https://arxiv.org/abs/2608.25175}{Optimization of closed-loop
  pumped-storage hydropower siting}, arXiv:2608.25175 [math.OC], accessed 31
  July 2026 (2026).
\newline\urlprefix\url{https://arxiv.org/abs/2608.25175}

\end{thebibliography}

\end{document}